\documentclass[a4paper,9pt]{article}
\usepackage{amsfonts}
\usepackage{amssymb,amsthm}
\usepackage{graphicx,amstext,amsmath,latexsym,dsfont,comment,color}

\newcommand{\R}{\mathbb{R}}
\newcommand{\N}{\mathbb{N}}

\newcommand{\Z}{\mathbb{Z}}

\newtheorem{theorem}{Theorem}

\newtheorem{example}{Example}

\newcommand{\be}{\begin{equation}}
\newcommand{\ee}{\end{equation}}
\newcommand{\bee}{\begin{equation*}}
\newcommand{\eee}{\end{equation*}}
\newcommand{\bp}{\begin{proof}}
\newcommand{\ep}{\end{proof}}

\date{20 June 2016 \\
\vspace {10mm}
{\it In memory of Yuri Safarov}}
\begin{document}
\title{On the number of Courant-sharp Dirichlet eigenvalues}
\author{M. van den Berg, K. Gittins
\thanks{MvdB acknowledges support by The Leverhulme Trust
through International Network Grant \emph{Laplacians, Random Walks, Bose Gas,
Quantum Spin Systems}. KG was supported by an EPSRC DTA. Both authors wish to thank Asma Hassannezhad for enjoyable discussions, and the referee for her/his helpful comments.}
 \\
\\
School of Mathematics, University of Bristol\\
University Walk, Bristol BS8 1TW, UK\\
\texttt{mamvdb@bristol.ac.uk}\\
\texttt{kg13951@bristol.ac.uk}
}

\maketitle

\begin{abstract}
We consider arbitrary open sets $\Omega$ in Euclidean space with finite Lebesgue measure, and obtain upper bounds for (i) the largest Courant-sharp Dirichlet eigenvalue of $\Omega$, (ii) the number of Courant-sharp Dirichlet eigenvalues of $\Omega$. This extends recent results of P. B\'erard and B. Helffer.
\end{abstract}

\section{Introduction}
\label{sec1}
Let $\Omega$ be an open set in Euclidean space  $\R^m$ with finite Lebesgue measure $|\Omega|$ and boundary $\partial \Omega$. We denote the spectrum of the Dirichlet Laplacian acting in $L^2(\Omega)$ by $\lambda_1(\Omega) \le \lambda_2(\Omega)\le \lambda_3(\Omega)\le \dots$ taking the multiplicities of these eigenvalues into account. We define the counting function for $\Omega$ by
\begin{equation*}
N_{\Omega}(\lambda)=\sharp\{n\in \N:\lambda_n(\Omega)<\lambda\}.
\end{equation*}
Weyl's law asserts that
\begin{equation}\label{e2}
N_{\Omega}(\lambda)=\frac{\omega_m}{(2\pi)^m}|\Omega|\lambda^{m/2}+o(\lambda^{m/2}),\, \lambda\rightarrow \infty,
\end{equation}
where $\omega_m$ is the measure of a ball $\mathcal{B}_m$ with radius $1$ in $\R^m$. We refer to Theorem 2 in \cite{GVR1} for a proof of \eqref{e2} in this generality.
For a proof of Weyl's law with a non-trivial remainder estimate for $\Omega$ open, bounded and connected we refer to Theorem 1.8 in \cite{NS}.

Let $\{\varphi_{1,\Omega},\varphi_{2,\Omega},\dots\}$ be an orthonormal basis in the Sobolev space $H_0^1(\Omega)$ of eigenfunctions corresponding to the Dirichlet eigenvalues. These eigenfunctions satisfy the Dirichlet boundary conditions in the usual trace sense. Let $\nu(\varphi_{n,\Omega})$ denote the number of nodal domains of $\varphi_{n,\Omega}$.
Then Pleijel's theorem (\cite{Pl}) states that
\begin{equation*}
\limsup_{n\rightarrow \infty} \frac{\nu(\varphi_{n,\Omega})}{n}\le \gamma_m,
\end{equation*}
where
\begin{equation}\label{e4}
\gamma_m=\frac{(2\pi)^m}{\omega_m^2}\big(\lambda_1(\mathcal{B}_m)\big)^{-m/2}<1.
\end{equation}
It is known that Pleijel's bound is not sharp. See \cite{B}, \cite{S} and \cite{P}.

We say that $\lambda_n(\Omega)$ is Courant-sharp if $\nu(\varphi_{n,\Omega})=n.$
Courant's nodal domain theorem asserts that $\nu(\varphi_{n,\Omega})\le n$. Courant's original proof in \cite{CH} was for the planar case. This has been subsequently stated and proved in a Riemannian manifold setting in \cite{C}. See also \cite{Pl}.
Pleijel's theorem implies that for a given $\Omega$ the number of Courant-sharp Dirichlet eigenvalues is finite. Using results of \cite{MvdBL} and \cite{Sa}, B\'erard and Helffer, \cite{BH}, obtained an upper bound for the largest Courant-sharp Dirichlet eigenvalue if $\Omega$ is bounded and has smooth boundary $\partial \Omega$.

This paper concerns arbitrary open sets in $\R^m$ with finite Lebesgue measure. The proofs of Courant's theorem in \cite{CH}, \cite{Pl} and \cite{C}) all use the fact that a restriction of an eigenfunction to a nodal domain $U$ is the first Dirichlet eigenfunction on $U$. This is immediate if $(\partial \Omega) \cap (\partial U)$ is sufficiently regular. The above fact holds without that regularity requirement. See for example Theorem 1.1 in \cite{H}.

Our main result, Theorem \ref{the1} below is for open sets $\Omega$ in $\R^m$ with finite Lebesgue measure.  We obtain (i) an upper bound for the largest Dirichlet eigenvalue of $\Omega$ which is Courant-sharp, and (ii) an upper bound for the number of Courant-sharp eigenvalues of $\Omega$. For $A\subset \R^m, A\ne \emptyset$ let
\begin{equation*}
d(x,A)=\inf\{|x-y|: y\in A\}.
\end{equation*}
For $\epsilon\ge 0$ and $|\Omega|<\infty$ we define
\begin{equation*}
\mu_{\Omega}(\epsilon)=|\{x\in \Omega: d(x,\partial \Omega)<\epsilon\}|,
\end{equation*}
and
\begin{equation}\label{e7}
\epsilon(\Omega)=\inf\{\epsilon:\mu_{\Omega}(\epsilon)\ge 2^{-1}(1-\gamma_m)|\Omega|\}.
\end{equation}
We denote the number of Courant-sharp eigenvalues of $\Omega$ by $\mathfrak{C}(\Omega)$.
\begin{theorem}\label{the1}
Let $\Omega$ be an open set in $\R^m$ with finite Lebesgue measure. We have the following.
\begin{enumerate}

\item[\textup{(i)}] If $\lambda_n(\Omega)$ is Courant-sharp then
\begin{equation}\label{e21}
\lambda_n(\Omega)\le \bigg(\frac{2\pi m^2}{(1-\gamma_m)\epsilon(\Omega)}\bigg)^2.
\end{equation}
\item[\textup{(ii)}]
\begin{equation}\label{e8}
\mathfrak{C}(\Omega)\le \frac{\omega_m}{(1-\gamma_m)^m}\big(m^3(m+2)\big)^{m/2}\frac{|\Omega|}{\epsilon(\Omega)^m}.
\end{equation}
\item[\textup{(iii)}] If $n\in \N, n>\frac{\omega_m}{(1-\gamma_m)^m}\big(m^3(m+2)\big)^{m/2}\frac{|\Omega|}{\epsilon(\Omega)^m},$ then $\lambda_n(\Omega)$ is not Courant-sharp.
\end{enumerate}
\end{theorem}

In Section \ref{sec2} below we prove Theorem \ref{the1}. In Section \ref{sec3} we analyse some examples including the von Koch snowflake.

\section{Proof of Theorem \ref{the1}}
\label{sec2}
Suppose $\lambda_n(\Omega)$ is Courant-sharp with eigenfunction $\varphi_{n,\Omega}$. Let $U_1,\dots,U_n$ be the nodal domains of $\varphi_{n,\Omega}$ so that $\lambda_n(\Omega)=\lambda_1(U_1)=\dots=\lambda_1(U_n)$. Without loss of generality
we may assume that $|U_1|\le |U_2|\le \dots \le |U_n|$. Hence  $|U_1|\le |\Omega|/n$. By Faber-Krahn we have that
\begin{equation*}
\lambda_n(\Omega)=\lambda_1(U_1)\ge \lambda_1(\mathcal{B}_m)\bigg(\frac{n\omega_m}{|\Omega|}\bigg)^{2/m}.
\end{equation*}
It follows that, since $\lambda_{n-1}(\Omega)<\lambda_{n}(\Omega)$,
\begin{align*}
(\lambda_n(\Omega))^{m/2}&\ge (\lambda_1(\mathcal{B}_m))^{m/2}\frac{n\omega_m}{|\Omega|}\nonumber \\
&\ge(\lambda_1(\mathcal{B}_m))^{m/2}\frac{\omega_m}{|\Omega|}(n-1)\nonumber \\ &
=(\lambda_1(\mathcal{B}_m))^{m/2}\frac{\omega_m}{|\Omega|}N_{\Omega}(\lambda_n(\Omega)).
\end{align*}
This gives that
\begin{equation}\label{e9}
\frac{\omega_m}{(2\pi)^m}(1-\gamma_m)|\Omega|(\lambda_n(\Omega))^{m/2}\le R_{\Omega}(\lambda_n(\Omega)),
\end{equation}
where $R_{\Omega}:\R^+\mapsto \R$ is defined by
\begin{equation}\label{e10}
R_{\Omega}(\lambda)=\frac{\omega_m}{(2\pi)^m}|\Omega|\lambda^{m/2}-N_{\Omega}(\lambda).
\end{equation}
See (15) and (16) in \cite{BH}. Below we obtain an upper bound for $R_{\Omega}(\lambda)$.
Let $\epsilon>0$ be arbitrary. Consider the collection $\mathfrak{M}_\epsilon$ of open cubes of measure $\epsilon^m$ with vertices in the set of $m$-tuples
$\{\Z\epsilon,\dots,\Z\epsilon\}$. Let $M_{\Omega}(\epsilon)$ be the number of open cubes of side-length $\epsilon$ in $\mathfrak{M}_{\epsilon}$ which are contained in $\Omega$,
\begin{equation*}
M_{\Omega}(\epsilon)=\sharp\{N\in \mathfrak{M}_{\epsilon}:N\subset \Omega\}.
\end{equation*}
We have that
\begin{equation}\label{e12}
|\Omega|-M_{\Omega}(\epsilon)\epsilon^m\ge 0.
\end{equation}
In order to obtain an upper bound for the left hand-side of \eqref{e12} we let $x\in \Omega$. If $d(x,\partial \Omega)>m^{1/2}\epsilon$, then $x$ belongs to an open $\epsilon$-cube in  $\mathfrak{M}_{\epsilon}$ contained in $\Omega$. Hence the measure of the set which is not covered by the $\epsilon$-cubes in $\mathfrak{M}_{\epsilon}$ that are entirely contained in $\Omega$ is bounded from above by $\mu_{\Omega}(m^{1/2}\epsilon)$. So
\begin{equation}\label{e13}
|\Omega|-M_{\Omega}(\epsilon)\epsilon^m\le\mu_{\Omega}(m^{1/2}\epsilon).
\end{equation}
By Dirichlet bracketing (see \cite{RS}) we have that
\begin{equation}\label{e14}
N_{\Omega}(\lambda)\ge M_{\Omega}(\epsilon)N_{C_\epsilon}(\lambda),
\end{equation}
where $C_{\epsilon}$ is an open cube in $\R^m$ with side-length $\epsilon$.
The following standard estimate is attributed to Gauss:
\begin{align}\label{e15}
N_{C_\epsilon}(\lambda)&=\sharp\bigg\{(k_1,\dots,k_m)\in \N^m: \sum_{i=1}^mk_i^2<\pi^{-2}\epsilon^2\lambda\bigg\}\nonumber \\
&\ge \frac{\omega_m}{2^m}\bigg(\pi^{-1}\epsilon\lambda^{1/2}-m^{1/2} \bigg)_+^m\nonumber \\
&\ge \frac{\omega_m}{(2\pi)^m}\epsilon^m\lambda^{m/2}\bigg(1-\frac{\pi m^{3/2}}{\epsilon \lambda^{1/2}} \bigg),
\end{align}
where $+$ denotes the positive part. By \eqref{e14} and \eqref{e15},
\begin{align}\label{e16}
N_{\Omega}(\lambda)&\ge M_{\Omega}(\epsilon)N_{C_\epsilon}(\lambda)\nonumber \\ &\ge
M_{\Omega}(\epsilon)\frac{\omega_m}{(2\pi)^m}\epsilon^m\lambda^{m/2}-M_{\Omega}(\epsilon)\frac{\omega_m}{(2\pi)^m}\pi m^{3/2}\epsilon^{m-1}\lambda^{(m-1)/2}\nonumber \\ &
=\frac{\omega_m}{(2\pi)^m}|\Omega|\lambda^{m/2}-(|\Omega|-M_{\Omega}(\epsilon)\epsilon^m)\frac{\omega_m}{(2\pi)^m}\lambda^{m/2}
\nonumber \\ &\hspace{25mm} -M_{\Omega}(\epsilon)\frac{\omega_m}{(2\pi)^m}\pi m^{3/2}\epsilon^{m-1}\lambda^{(m-1)/2}.
\end{align}
We bound the second and third terms in the right hand-side of \eqref{e16} using \eqref{e13} and \eqref{e12} respectively. This then gives, by \eqref{e10}, that
\begin{equation}\label{e17}
R_{\Omega}(\lambda)\le \frac{\omega_m}{(2\pi)^m}\mu_{\Omega}(m^{1/2}\epsilon)\lambda^{m/2}+\frac{\pi m^{3/2}\omega_m}{(2\pi)^m}\frac{|\Omega|\lambda^{(m-1)/2}}{\epsilon}.
\end{equation}
By \eqref{e9} and \eqref{e17} we have that if $\lambda_n(\Omega)$ is Courant-sharp then
\begin{align}\label{e18}
\frac{\omega_m}{(2\pi)^m}(1-\gamma_m)|\Omega|(\lambda_n(\Omega))^{m/2}&\le  \frac{\omega_m}{(2\pi)^m}\mu_{\Omega}(m^{1/2}\epsilon)(\lambda_n(\Omega))^{m/2}\nonumber\\ & \hspace{5mm}+\frac{\pi m^{3/2}\omega_m}{(2\pi)^m}\frac{|\Omega|(\lambda_n(\Omega))^{(m-1)/2}}{\epsilon}.
\end{align}
We now choose $\epsilon$ such that the second term in the right hand-side of \eqref{e18} equals half of the left hand-side of \eqref{e18}. That is
\begin{equation}\label{e19}
\epsilon=2\pi m^{3/2}(1-\gamma_m)^{-1}(\lambda_n(\Omega))^{-1/2}.
\end{equation}
By \eqref{e18} and the choice of $\epsilon$ in \eqref{e19} we have that if $\lambda_n(\Omega)$ is Courant-sharp then
\begin{equation}\label{e20}
2^{-1}(1-\gamma_m)|\Omega|\le \mu_{\Omega}(2\pi m^2(1-\gamma_m)^{-1}(\lambda_n(\Omega))^{-1/2}).
\end{equation}
Since $\epsilon\mapsto \mu_{\Omega}(\epsilon)$ is continuous and onto $[0,|\Omega|]$ the infimum in \eqref{e7} is over a non-empty set which is bounded from below, and therefore exists.
So if $\lambda_n(\Omega)$ is Courant-sharp then, by \eqref{e7} and \eqref{e20}, $\frac{2\pi m^2}{(1-\gamma_m)(\lambda_n(\Omega))^{1/2}}\ge \epsilon(\Omega)$.
This proves Theorem \ref{the1}(i).

By \cite{LY} we also have that
\begin{equation}\label{e22}
\lambda_n(\Omega)\ge \frac{m}{m+2}\frac{(2\pi)^2}{\omega_m^{2/m}}\bigg(\frac{n}{|\Omega|}\bigg)^{2/m}.
\end{equation}
This, together with \eqref{e21}, implies \eqref{e8} and proves Theorem \ref{the1}(ii).

To prove Theorem \ref{the1}(iii) we just note that by \eqref{e22}, $$\max\bigg\{n\in \N:\lambda_n(\Omega)\le \bigg(\frac{2\pi m^2}{(1-\gamma_m)\epsilon(\Omega)}\bigg)^2\bigg\}\le \frac{\omega_m}{(1-\gamma_m)^m}\big(m^3(m+2)\big)^{m/2}\frac{|\Omega|}{\epsilon(\Omega)^m}.$$

\hspace*{\fill }$\square $

We note that if we were to use the lower bounds for the counting function from Section 2 in \cite{MvdBL} then we would have to assume
a weak integrability condition on $\mu_{\Omega}$ of the form $\int \epsilon^{-1}\,d\mu_{\Omega}(\epsilon)<\infty$. Such an integrability condition may fail if the interior Minkowski dimension of $\partial \Omega$ is equal to $m$. The procedure above avoids this integrability condition.

\section{Examples}
\label{sec3}

In this section we analyse three examples where explicit computations seem out of reach.

\begin{example}\label{ex1}
Let $\Omega$ be an open, bounded, convex set in $\R^m$. Let $\mathcal{H}^{m-1}(\partial \Omega)$ denote the $(m-1)$-dimensional Hausdorff measure of $\partial \Omega$.
Then \begin{equation}\label{e41}
\mathfrak{C}(\Omega)\le \frac{\omega_m}{(1-\gamma_m)^{2m}}\big(4m^3(m+2)\big)^{m/2}\frac{\big(\mathcal{H}^{m-1}(\partial \Omega)\big)^m}{|\Omega|^{m-1}}.
\end{equation}
\begin{proof}
By convexity of $\Omega$ we have that
\begin{equation*}
\mu_{\Omega}(\epsilon)\le \mathcal{H}^{m-1}(\partial \Omega)\epsilon.
\end{equation*}
By \eqref{e7},
\begin{equation}\label{e43}
\epsilon(\Omega)\ge 2^{-1}(1-\gamma_m)\frac{|\Omega|}{\mathcal{H}^{m-1}(\partial \Omega)},
\end{equation}
and \eqref{e41} follows from Theorem \ref{the1} and \eqref{e43}.
\end{proof}
\end{example}

It was shown in \cite{HS} that only the first, second and fourth Dirichlet eigenvalues for $\mathcal{B}_2$ are Courant-sharp. Hence $\mathfrak{C}(\mathcal{B}_2)=3,$ and the largest Courant-sharp eigenvalue for $\mathcal{B}_2$ is equal to $j_{0,2}^2$. Here $j_{0,2}\asymp 5.520..$ is the second positive zero of the Bessel function $J_0$. A straightforward computation using \eqref{e21} and \eqref{e43} shows that the largest Courant-sharp eigenvalue of $\mathcal{B}_2$  is strictly less than $1.2\cdot10^6$. This compares well with the bound $7.1\cdot 10^6$ obtained in \cite{BH}. For the unit square $\mathcal{C}_2$ it is known (\cite{Pl}, \cite{BH1}) that only the first, second and fourth Dirichlet eigenvalues are Courant-sharp. Hence $\mathfrak{C}(\mathcal{C}_2)=3,$
and the largest Courant-sharp eigenvalue for $\mathcal{C}_2$ is equal to $8\pi^2$. Using \eqref{e21} and \eqref{e43} we have that the largest Courant-sharp eigenvalue of the unit square is strictly less than $4.5\cdot10^6,$ whereas \cite{BH} gives a bound $5.9\cdot10^6.$ These examples illustrate that the bounds obtained in Theorem \ref{the1} are very crude. 

\begin{figure}[h]
\centering\includegraphics[scale=.7]{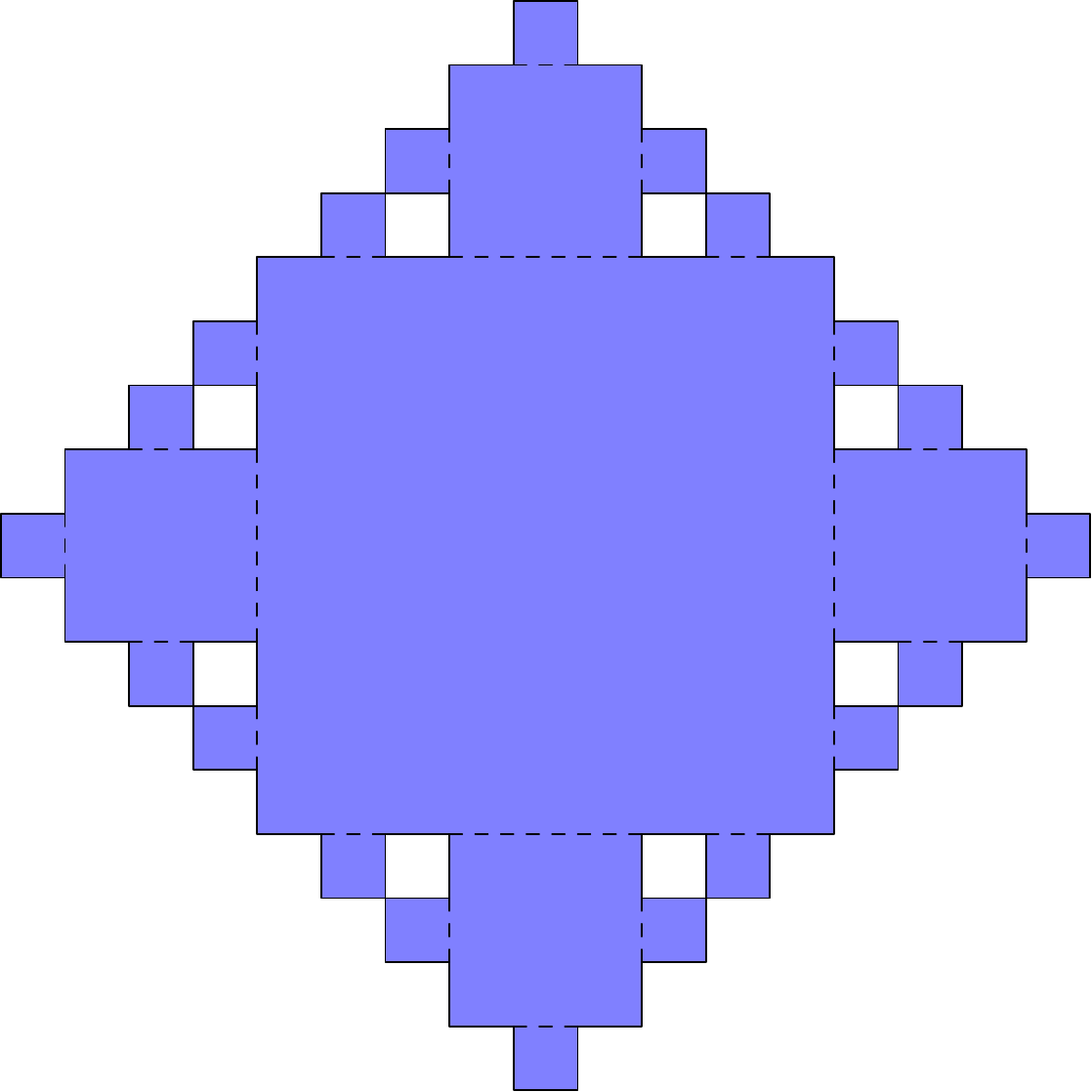}
\caption{The first two generations of $K$}\label{fig1}
\end{figure}

The second example is a von Koch snowflake $K$ with similarity ratio $\frac13$. We recall its construction. Let the basic square (generation $0$) in $K$ have side-length $1$. The first generation consists of $4$ squares with side-length $\frac13$ each attached symmetrically to the basic square.
Proceeding inductively we have that the $j$'th generation in $K$, $j\in \N$ consists of $4\cdot5^{j-1}$ squares with side-length $3^{-j}$. We let $K$ be the interior of its closure. Then $K$ is connected, has Lebesgue measure $|K|=2$, and both the Hausdorff dimension of $\partial K$ and the interior Minkowski dimension of $\partial K$ are equal to $\log 5/ \log 3.$ See Figure~\ref{fig1}, and  \cite{MvdB} for further details.

\begin{example}\label{ex2} Let $K$ be the von Koch snowflake generated by the unit square and similarity ratio $\frac13$. Then
\begin{equation}\label{e23}
\mathfrak{C}(K)\le 15\cdot10^7.
\end{equation}
\end{example}

\begin{proof}
By Theorem \ref{the1}, \eqref{e4}, and $|K|=2$ we find that
\begin{equation}\label{e24}
\mathfrak{C}(K)\le\frac{64\pi j_0^4}{(j_0^2-4)^2}\epsilon(K)^{-2},
\end{equation}
where we have used that
\begin{equation*}
\lambda_1(\mathcal{B}_2)=j_0^2,
\end{equation*}
where $j_0=2.405...$ is the first positive zero of the Bessel function $J_0$. It remains to find a lower bound for $\epsilon(K)$. We obtain an upper bound for $\mu_{\Omega}(\epsilon)$ by adding all edges between squares of different generations. This gives a disjoint union of $1$ unit square and $4\cdot5^{j-1}$ squares with side-lengths $3^{-j},j\in \N$. Let $\epsilon<\frac{1}{18}$, and let $J\in \N$ be such that
\begin{equation*}
J< \frac{\log\big(\frac{1}{2\epsilon}\big)}{\log 3}\le J+1.
\end{equation*}
Then $J\ge 2$. The contribution to the upper bound for $\mu_{\Omega}(\epsilon)$ from the squares in generations $1,\dots,J-1$ is bounded from above by
\begin{equation}\label{e27}
\bigg(4+ 16 \sum_{j=1}^{J-1}5^{j-1}3^{-j}\bigg)\epsilon\le \frac{24\epsilon}{5}\bigg(\frac{5}{3}\bigg)^J\le \frac{48}{5}2^{-\frac{\log 5}{\log 3}} \epsilon^{2-\frac{\log 5}{\log 3}}.
\end{equation}
The first term in the left-hand side above is the contribution from the unit square. The contribution to the upper bound for $\mu_{\Omega}(\epsilon)$ from the squares in generations $J,J+1,\dots$ is bounded from above by
\begin{equation}\label{e28}
\sum_{j\ge J} 4\cdot5^{j-1}9^{-j}=\bigg(\frac{5}{9}\bigg)^{J-1}\le \frac{36}{5}2^{-\frac{\log 5}{\log 3}} \epsilon^{2-\frac{\log 5}{\log 3}}.
\end{equation}
We recognise the interior Minkowski dimension $\frac{\log 5}{\log 3}$ of $\partial K$. By \eqref{e27} and \eqref{e28} we have that
\begin{equation*}
\mu_{\Omega}(\epsilon)\le \frac{84}{5}2^{-\frac{\log 5}{\log 3}} \epsilon^{2-\frac{\log 5}{\log 3}},\, 0<\epsilon<\frac{1}{18}.
\end{equation*}
Solving the equation
\begin{equation*}
 \frac{84}{5}2^{-\frac{\log 5}{\log 3}} \epsilon^{2-\frac{\log 5}{\log 3}}= 1-\frac{4}{j_0^2}
\end{equation*}
gives that
\begin{equation}\label{e31}
\epsilon(K)\ge 0.00379.
\end{equation}
The bound of \eqref{e23} follows by \eqref{e24} and \eqref{e31}.
\end{proof}

\begin{figure}[h]
\centering\includegraphics[scale=1.0]{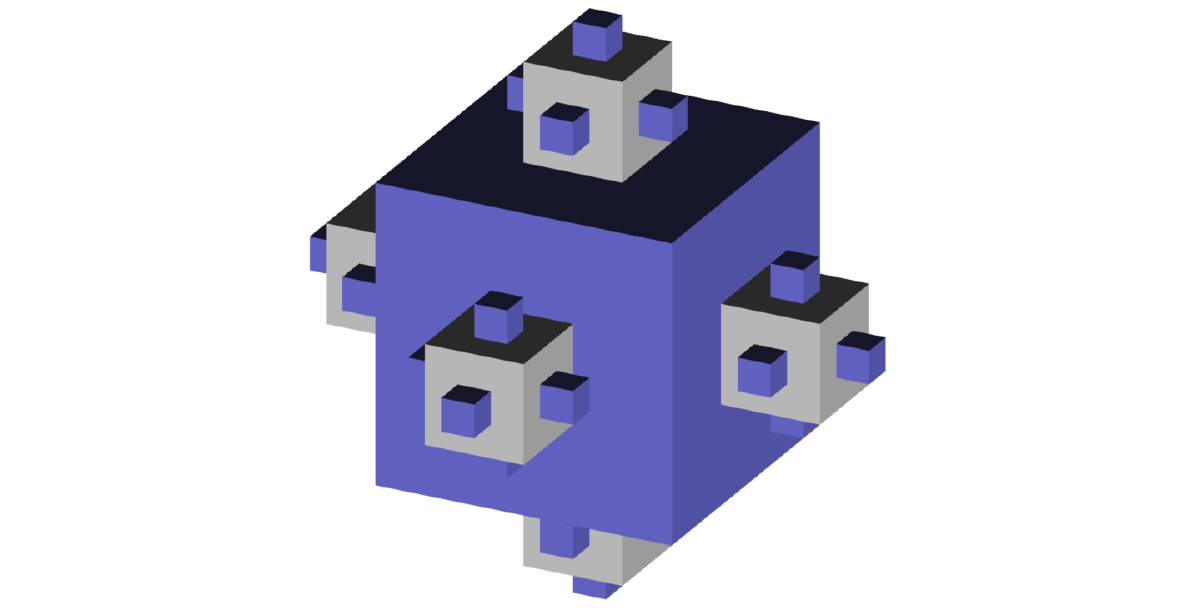}
\caption{The first two generations of $D_{s}$ with $s=\frac{1}{3}$.}\label{fig2}
\end{figure}

Below we construct an open set $D_s\subset \R^3$. Let $Q_{0} \subset \R^{3}$ be an open cube of
side-length 1. Let $0<s\le \sqrt 2-1$. Attach a regular open cube $Q_{1,i}$ of side-length $s$ to the
centre $c_{1,i}, i=1,\dots,6,$ of each face of $\partial Q_{0}$, and such that all the faces are pairwise-parallel.
Now proceed by induction. For $j=2,3,\dots,$ attach $N(j)=6\cdot5^{j-1}$ open cubes $Q_{j,1}, \dots, Q_{j,N(j)},$
of side-length $s^{j}$ to the centres of the boundary faces of the cubes $Q_{j-1,1},\dots,Q_{j-1,N(j-1)}$,
again with pairwise-parallel faces. We define the polyhedron $D_{s}$ as
\begin{equation*}
D_{s}=\text{interior}\bigg\{\overline{ Q_{0} \cup [\bigcup_{j\geq1} \bigcup_{1\leq i \leq N(j)} Q_{j,i} ]}\bigg\}.
\end{equation*}
See Figure~\ref{fig2}. We note that for $0<s\le\sqrt 2-1$ no cubes in the construction of $D_{s}$ overlap.

The asymptotic behaviour of the heat content of $D_s$ in $\R^3$ for small time was analysed in \cite{vdBG}. Here we have the following.
\begin{example}\label{ex3} Let $s\in (0,\sqrt 2-1]$, and let $D_s$ be the polyhedron in $\R^3$ defined above. Then
\begin{equation}\label{e33}
\mathfrak{C}(D_s)\le 25\cdot10^{10}.
\end{equation}
\end{example}
\begin{proof}
We have that
\begin{equation*}
\vert D_{s} \vert = \frac{1+s^{3}}{1-5s^{3}}, \
\end{equation*}
and that the two-dimensional Hausdorff measure of the boundary is given by
\begin{equation*}
\mathcal{H}^2(\partial D_{s})  = 6\left(\frac{1-s^{2}}{1-5s^{2}}\right).
\end{equation*}
By Theorem \ref{the1} we have that
\begin{equation}\label{e36}
\mathfrak{C}(D_s)\le 36(15)^{3/2}\pi\bigg(1-\frac{9}{2\pi^2}\bigg)^{-3}\frac{|D_s|}{\epsilon(D_s)^3},
\end{equation}
where we have used that
\begin{equation*}
\lambda_1(\mathcal{B}_3)=j_{1/2}^2=\pi^2,
\end{equation*}
where $j_{1/2}=\pi$ is the first positive zero of the Bessel function $J_{1/2}$.
We obtain an upper bound for $\mu_{\Omega}(\epsilon)$ by adding all faces between cubes of different generations.
This gives a disjoint union of $1$ unit cube and $6\cdot5^{j-1}$ cubes of side-length $s^j,\, j\in \N$. Hence
\begin{equation}\label{e38}
\mu_{\Omega}(\epsilon)\le \bigg(6+36\sum_{j=1}^\infty5^{j-1}s^{2j}\bigg)\epsilon=\frac{6(1+s^2)}{1-5s^2}\epsilon.
\end{equation}
By \eqref{e7} and \eqref{e38} we have that
\begin{equation}\label{e39}
\epsilon(D_s)\ge\frac{1}{12}\bigg(1-\frac{9}{2\pi^2}\bigg)\frac{1-5s^2}{1+s^2}|D_s|.
\end{equation}
Finally by \eqref{e36}, \eqref{e39}, the fact that $0<s\le \sqrt 2-1$, and $|D_s|\ge 1$ we obtain that
\begin{equation*}
\mathfrak{C}(D_s)\le 6(12)^4(15)^{3/2}(140+99\sqrt 2)\pi\bigg(1-\frac{9}{2\pi^2}\bigg)^{-6}.
\end{equation*}
This implies \eqref{e33}.
\end{proof}

\end{document}